\documentclass[a4paper,12pt]{article}
\title{ A Simple Proof of the Classification of Normal Toeplitz Matrices}
\author{Akio Arimto\footnote{ Department of Electronics and Computer Engineering, Musashi Institute of Technology, Tamazutsumi 1-28-1,Setagaya-ku, Tokyo 158,Japan(arimoto@iname.com)}}
\date{}
\begin{document}
\maketitle
\begin{abstract} 
 We give an easy proof to show that every complex normal Toeplitz matrix is classified as either of type I or of type II. 
 Instead of difference equations on elements in the matrix used in past studies, polynomial equations with coefficients of elements are used. 
 In a similar fashion, we show that a real normal Toeplitz matrix must be one of four types: symmetric, skew-symmetric, circulant, or skew-circulant. 
 Here we use trigonometric polynomials in the complex case and algebraic polynomials in the real case.
\end{abstract}

\section
{Introduction. Normal Toeplitz Matrices. }

Ikramov[3] has firstly shown that a normal Toeplitz matrix over real field must be one of four types: symmetric, skew-symmetric (up to the principal diagonal), circulant, or skew-circulant. 
Later on the complex version of the problem has been solved by Ikramov[6] et.al and Gel'fgat[2]. 
The history of the theory has been described in [5]. 
The complex case was also proved by Farenick, Krupnik, Krupnik and Lee[1], and Ito[2] independently to [6] . 
We will give here another  proof which, the author believes, is the simplest one. 
The complex case and the real case will be proved in a similar fashion, although trigonometric polynomials will be used in the complex case and algebraic polynomials will be used in the real case, respectively.  
Now denote Toeplitz matrix of order $N + 1$ by $T_N $:
$$\,\,\,\,\,\,\,\,\,\,\,T_N \, = \left[ {\begin{array}{*{20}c}
   {a_0 } & {a_{ - 1} } &  {a_{-2}}&\cdots  & {a_{ - N} }  \\
   {a_1 } & {a_0 } & {a_{ - 1} } & \cdots  & {a_{ - \left( {N - 1} \right)} }  \\    \cdots  &  \cdots  &  \cdots  &  \cdots   \\
   {a_{N - 1} } & {a_{N - 2} } & {a_{N - 3} } &  \cdots  & {a_{ - 1} }  \\
   {a_N } & {a_{N - 1} } & {a_{N - 2} } &  \cdots  & {a_0 }  \\
\end{array}} \right].$$
$T_N $ is called normal Toeplitz if $T_N T_N^*  - T_N^* T_N  = 0,$ where $T_N^* $ is a transposed conjugate of $T_N $. 
We can simply state classifications of normal Toeplitz matrices in the following way. 
A normal Toeplitz matrix $T_N $ is of type I if and only if for some $\alpha _0 $,with $\left| {\alpha _0 } \right| = 1$,
$$\,\,\,\left[ {\begin{array}{*{20}c}
   {a_{ - 1} } & {a_{ - 2} } &  \cdots  & {a_{ - N} }  \\
\end{array}} \right] = \,\alpha _0
 \left[
  {\begin{array}{*{20}c}
   {\bar a_1 } & {\bar a_2 } &  \cdots  & {\bar a_N }  \\
\end{array}} \right],
 \leqno{(1)}$$
whereas $T_N$ is of type II if and only if for some $\beta _0$, 
with $\left|{\beta_0}\right|=1,$
$$\,\,\,\left[ {\begin{array}{*{20}c}
   {a_{ - 1} } & {a_{ - 2} } &  \cdots  & {a_{ - N} }  \\
\end{array}} \right] = \,\beta _0 \left[ {\begin{array}{*{20}c}
   {a_N } & {a_{N - 1} } &  \cdots  & {a_1 }  \\
\end{array}} \right].\leqno{(2)}$$
Surely the complex case contains the real case, but when all the $a_k 's$ are real valued, it will be found that $\alpha _0 $ and  $\beta_0$ must be  +1 or -1, hence we may have the further concrete classification of four types : A normal Toeplitz matrix  $T_N$ is symmetric if and only if
 $$\,\,\,\left[ {\begin{array}{*{20}c} {a_{ - 1} } & {a_{ - 2} } &  \cdots  & {a_{ - N} }  \\
\end{array}} \right] = \,\left[ {\begin{array}{*{20}c}
   {a_1 } & {a_2 } &  \cdots  & {a_N }  \\
\end{array}} \right],\leqno{(3)}$$
$T_N$  is skew-symmetric if and only if 
  $$\,\,\,\left[ {\begin{array}{*{20}c}
   {a_{ - 1} } & {a_{ - 2} } &  \cdots  & {a_{ - N} }  \\
\end{array}} \right] = -\,\left[ {\begin{array}{*{20}c}
   {a_1 } & {a_2 } &  \cdots  & {a_N }  \\
\end{array}} \right],\leqno{(4)}$$
$T_N$  is circulant if and only if
$$\,\,\,\left[ 
{\begin{array}{*{20}c}
   {a_{ - 1} } & {a_{ - 2} } &  \cdots  & {a_{ - N} }  \\
\end{array}}
 \right] = \,\left[
  {\begin{array}{*{20}c}
   {a_N } & {a_{N - 1} } &  \cdots  & {a_1 }  \\
 \end{array}}
  \right],\leqno{(5)}$$
$T_N$  is skew-circulant if and only if
 $$\,\,\,\left[ {\begin{array}{*{20}c}
   {a_{ - 1} } & {a_{ - 2} } &  \cdots  & {a_{ - N} }  \\
\end{array}} \right] = -\,\left[ {\begin{array}{*{20}c}
  {a_N } & {a_{N - 1} } &  \cdots  & {a_1 }  \\
\end{array}} \right]. \leqno{(6)}$$

Throughout this paper, we assume  $a_0$ to be zero without loss of generality since  $a_0$ does not play any role in the above classifications.

\section{Proof of Theorem in the Complex Case}

 Now we define trigonometric polynomials $s(x)$ and $t(x)$ with coefficients $a_1 ,a_2 , \cdots ,a_N $ and $a_{ - 1} ,a_{ - 2} , \cdots ,a_{ - N} $, respectively:
$$s\left( x \right) = a_1 e^{ix}  + a_2 e^{2ix}  +  \cdots  + a_N e^{Nix} $$
$$t\left( x \right) = a_{ - 1} e^{ - ix}  + a_{ - 2} e^{ - 2ix}  +  \cdots  + a_{ - N} e^{ - Nix} $$
Using these trigonometric polynomials, we can restate the conditions (1) and 
(2) in the previous section as \\
(1') type I: $t\left( x \right) = \alpha _0 \overline {s\left( x \right)} $
, with $\left| {\alpha _0 } \right|\, = 1$, \\
(2') type II: $t\left( x \right) = \beta _0 s\left( x \right)e^{ - i\left( {N + 1} \right)x} $, with $\left| {\beta _0 } \right| = 1$.
\\\\
\textbf{Theorem 1} \textit{Every normal Toeplitz matrix is either of type I or of type II.} \\\\
\textit{Proof.} Our proof is based on the expression in [7] p.998 or [1] 
p.1038. It is the necessary and sufficient conditions for $T_N$ to be normal in terms of $a_n$ that 
$$a_m \bar a_n  - \bar a_{ - m} a_{ - n}  + \bar a_{N + 1 - m} a_{N + 1 - n}  - a_{ - \left( {N + 1 - m} \right)} \bar a_{ - \left( {N + 1 - n} \right)}  = 0,\leqno{(7)}$$
$\left( {1 \le n,m \le N} \right). $ \\

We now rewrite (7) by using polynomials $s(x)$ and $t(x)$. 
Multiply both sides of (7) by $e^{imx} e^{ - iny} $ and sum these new expressions over all $m$ and $n$ to obtain
 $$s\left( x \right)\overline {s\left( y \right)} \, - \overline {t\left( x \right)} \,t\left( y \right) + \overline {s\left( x \right)\,} s\left( y \right)e^{i\left( {N + 1} \right)\left( {x - y} \right)}  - t\left( x \right)\overline {t\left( y \right)} e^{i\left( {N + 1} \right)\left( {x - y} \right)}  = 0
.  \leqno{(8)}
$$
It easily seen that (8) is equivalent to (7) , and (1') or (2') implies (8). 
Hence both (1') and (2') are sufficient conditions for $T_N$ to be normal. So it only remains to show the necessity of the condition.
When we take $x=y$ in (8), we have 

$$s\left( x \right)\overline {s\left( x \right)}  - t\left( x \right)\overline {t\left( x \right)}  = 0, or \left| {s\left( x \right)} \right| = \left| {t\left( x \right)} \right|  \leqno{(9)}$$
for each $x$. Hence except for trivial case $t\left( x \right) \equiv 0$
, we can find an $x_0$ such that $t\left( {x_0 } \right) \ne 0$
and there exist $\alpha $ and $\beta $ such that $s\left( {x_0 } \right) = \alpha t\left( {x_0 } \right)$
and $t\left( {x_0 } \right) = \beta \overline {t\left( {x_0 } \right)} $
, with $\left| \alpha  \right| = \left| \beta  \right| = 1$. By setting $y = x_0 $ and then by dividing through by $\overline {t\left( {x_0 } \right)}$, equation (8) becomes

$$s\left( x \right)\bar \alpha  - \overline {t\left( x \right)} \beta  + \overline {s\left( x \right)} \alpha \,\beta e^{i\left( {N + 1} \right)\left( {x - x_0 } \right)}  - t\left( x \right)e^{i\left( {N + 1} \right)\left( {x - x_0 } \right)}  = 0 \leqno{(10)}$$
from which we have

$$\overline {t\left( x \right)}  = s\left( x \right)\overline {\alpha \,\beta }  + \overline {s\left( x \right)} \alpha e^{i\left( {N + 1} \right)\left( {x - x_0 } \right)}  - t\left( x \right)\bar \beta e^{i\left( {N + 1} \right)\left( {x - x_0 } \right)}. \leqno{(11)}$$

Substituting the right hand side of (11) into $\overline {t\left( x \right)} $ of (9), we have
$$\left\{ {s\left( x \right) - t\left( x \right)\alpha e^{i\left( {N + 1} \right)\left( {x - x_0 } \right)} } \right\}\left\{ {\overline {s\left( x \right)}  - t\left( x \right)\overline {\alpha \beta } } \right\} = 0, \leqno{(12)}$$
which implies our required results if we take $\alpha _0  = \alpha \beta $
and $\beta _0  = \bar \alpha e^{i\left( {N + 1} \right)x_0 } $.
We should notice here that (12) is the product of polynomials. 
If one of the polynomials in (12) is nonzero at $x = x_0 $, then by continuity the polynomial is nonzero in an open neibourhood $\mathcal{U}$ of $x_0$, which implies that the other polynomial in the product (12) is identically zero on $\mathcal{U}$. 
Because the only polynomial with a continuum of roots is the zero polynomial, equation (12) implies equation (1') or (2'). 
Hence type I or type II must occur. 

\section{Proof of Theorem in the Real Case. }

\textbf{Theorem 2}\textit{ Every real normal Toeplitz matrix is either symmetric, skew-symetric, circulant or skew-circulant.}\\\\
\textit{Proof.}  We define algebraic polynomials $p(x)$ and $q(x)$ such that
$$p\left( x \right) = a_1 x + a_2 x^2  +  \cdots  + a_N x^N $$
$$q\left( x \right) = a_{-1} x + a_{-2} x^2  +  \cdots  + a_{-N} x^N $$
Also in this real case we can apply (7) in the previous section as necessary and sufficient conditions for $T_N$ to be normal in terms of real valued
$a_n$:
$$a_m a_n  - a_{ - m} a_{ - n}  + a_{N + 1 - m} a_{N + 1 - n}  - a_{ - \left( {N + 1 - m} \right)} a_{ - \left( {N + 1 - n} \right)}  = 0 \left( {1 \le n,m \le N} \right)
 \leqno{(13)}$$
 Now we introduce the reciprocal polynomials to $p(x)$ and $q(x)$ 

$$\tilde p\left( x \right) = a_N x + a_{N - 1} x^2  +  \cdots  + a_1 x^N $$
and
$$\tilde q\left( x \right) = a_{ - N} x + a_{ - N + 1} x^2  +  \cdots  + a_{ - 1} x^N, $$
respectively. Multiply both sides of equation (13) by $x^m y^n $ and then sum the resulting expressions over all $m$ and $n$ to obtain the following analogue of (8):
$$p\left( x \right)p\left( y \right) - q\left( x \right)q\left( y \right) + \tilde p\left( x \right)\tilde p\left( y \right) - \tilde q\left( x \right)\tilde q\left( y \right) = 0. \leqno{(14)} $$
Noticing the relation $\tilde q\left( x \right) = x^{N+1} q\left( {1/x} \right)$
 and $\tilde p\left( x \right) = x^{N+1} p\left( {1/x} \right)$
, and letting $y = 1/x$
  in (14), we have $p\left( x \right)\tilde p\left( x \right) = q\left( x \right)\tilde q\left( x \right)$ .
 By setting $x = y$ and then multiplying both sides of (14) by $p^2
   \left( x \right)$, we have from the relation $p^2 \left( x \right)\tilde p^2 \left( x \right) = q^2 \left( x \right)\tilde q^2 \left( x \right)$
 ,
$$p^4 \left( x \right) - p^2 \left( x \right)q^2 \left( x \right) + q^2 \left( x \right)\tilde q^2 \left( x \right) - p^2 \left( x \right)\tilde q^2 \left( x \right) = 0, \leqno{(15)}$$ 
which can be factorized as the following way:
$$\left( {p^2 \left( x \right) - q^2 \left( x \right)} \right)\left( {p^2 \left( x \right) - \tilde q^2 \left( x \right)} \right) = 0. \leqno{(16)}$$
From this formula, we obtain four types, symmetric $p\left( x \right) = q\left( x \right)$
 , skew-symmetric $p\left( x \right) =  - q\left( x \right)$
 , circulant $p\left( x \right) = \tilde q\left( x \right)$
  and skew-circulant $p\left( x \right) =  - \tilde q\left( x \right)$
 , which are equivalent to (3),(4),(5) and (6), respectively.

\end{document}